\documentclass[11pt,twoside,letter]{article}
\usepackage{graphicx}
\usepackage{fancyhdr}
\usepackage{amsmath}
\usepackage{amssymb}
\usepackage{enumerate}

\newcommand\bk{{\mathbf k}}
\newcommand\e{{\varepsilon}}
\newcommand\cV{{\mathcal V}}
\newcommand\R{{\mathbf R}}
\newcommand\N{{\mathbf N}}
\newcommand\Z{{\mathbf Z}}
\newcommand\cF{{\mathcal F}}
\newcommand{\be}{\begin{equation}}
\newcommand{\ee}{\end{equation}}
\newcommand{\iref}[1]{{\rm (\ref{#1})}}
\newcommand\eref[1]{{\rm (\ref{#1})}}

\DeclareMathOperator*{\argmin}{argmin}

\newtheorem{lemma}{Lemma}[section]
\newtheorem{theorem}[lemma]{Theorem}

\begin{document}
\title
{
Sampling
and reconstruction of solutions to the Helmholtz equation
}
\author{ 
Gilles Chardon \\ \small Acoustics Research Institute, Austrian Academy of Sciences\\ \small A-1040, Wien, Austria \\ \small gilles.chardon@m4x.org\\
Albert Cohen \\ \small UPMC Univ. Paris 06, UMR 7598, Laboratoire Jacques-Louis Lions\\ \small F-75005, Paris, France \\  \small albert.cohen@upmc.fr\\
Laurent Daudet \\ \small Institut Langevin, Paris Diderot University and Institut Universitaire de France\\ \small F-75005, Paris, France \\ \small laurent.daudet@espci.fr
}

\date{}

\maketitle
\thispagestyle{fancy}

\markboth{\footnotesize \rm \hfill G. CHARDON, A. COHEN AND  L. DAUDET}
{\footnotesize \rm \hfill Sampling
and reconstruction of solutions to the Helmholtz equation \hfill}

\begin{abstract}
We consider the inverse problem of reconstructing 
general solutions
to the Helmholtz equation  on some domain $\Omega$
from their values at scattered points
$x_1,\dots,x_n\subset \Omega$.
This problem typically arises when sampling acoustic 
fields with $n$ microphones for the purpose of
reconstructing this field over a region of interest $\Omega$ 
contained in a larger domain $D$ in which the acoustic field propagates.
In many applied settings, the shape of $D$ and the boundary conditions
on its border are unknown.
Our reconstruction method is based on
the approximation of a general solution $u$
by linear combinations of Fourier-Bessel functions or plane waves.
We analyze the convergence of 
the least-squares estimates to $u$ 
using these families of functions based
on the samples $(u(x_i))_{i=1,\dots,n}$. Our analysis describes the 
amount of regularization needed to 
guarantee the convergence of the least squares estimate
towards $u$, in terms of a condition that depends on the
dimension of the approximation subspace,
the sample size $n$ and
the distribution of the samples. It reveals the
advantage of using non-uniform distributions
that have more points on the boundary of $\Omega$.
Numerical illustrations show that our approach compares
favorably with reconstruction methods
using other basis functions, and other types
of regularization.
\vspace{5mm}\\
\noindent
 {\it Key words and phrases} : Helmholtz equation, interpolation, least squares, regularization
\vspace{3mm}\\
\noindent
MSC2000: 74J25,35J05,94A20
\end{abstract}

\section{Introduction}  

A common inverse problem in acoustics is to obtain a precise approximation 
of the soundfield over a spatial domain $\Omega$ of interest,
using the smallest possible number 
of pointwise measurements, e.g.\ as provided by microphones. 
For instance, one may wish to measure the complex radiation pattern of 
an extended source (source identification problem), to localize a number of point sources within a spatial 
domain (source localization problem), or to optimize the output of a sound reproduction system over a large control area, to name only a few applications. 
In practice, the main difficulty that one is usually faced with is how to handle reverberation: the reverberant field might well be of a magnitude comparable to the direct sound, and it depends in a non-trivial way on both the geometry of the domain $D$ where the acoustic field is defined 
and the type of boundary conditions 
on $\partial D$
(with Dirichlet or Neumann as ideal cases, but more likely in engineering problems
 with a frequency-dependent mixed behavior). 

The goal of this paper is to study the accuracy that can be achieved
when approximating the acoustic field over the domain $\Omega\subset D$,
based on a set of point measurements, 
{\it without precise knowledge on the geometry of $D$ and boundary conditions on $\partial D$}. 
A general setting is the following: the soundfield $p(x,t)$ is measured at microphones located
at positions $x_1,\dots,x_n\in \Omega$,
and over a (discretized) time interval $[0,T]$. After application of the (discrete) Fourier transform $\cF$
in the time variable, and considering a given frequency $\omega$, the function
$$
u(x):=\cF p(x,\omega)
$$
is a solution on $D$
to the Helmholtz equation 
\be
\label{helm}
\Delta u+\lambda^2 u=0,
\ee
where $\lambda = \omega / c$, with $c$ denoting the wave velocity,
and where the boundary conditions are unknown to us. 

Depending on the applications, the geometry of the domain $\Omega$ may either be
$2$-D (membranes) or $3$-D (rooms). Our problem therefore amounts to reconstructing, on some domain $\Omega\subset \R^2 \mbox{ or } \R^3$,
a general solution to the Helmholtz equation \iref{helm}
from its sampling at points $x_1,\dots,x_n\in \Omega$. These
samples may be measured exactly or up to some additive noise. We denote by
\be
y_l=u(x_l)+\eta_l,\quad\quad l=1,\dots,n,
\label{data}
\ee
these samples, where $\eta_l$ represent the additive noise.

Reconstruction from scattered points is a widely studied topic,
and a variety of methods have been proposed and analyzed.
Many existing methods can be viewed as reconstructing
some form of approximation to the unknown function $u$
by simpler functions such as splines, partial Fourier sums
or radial basis functions. The success of these methods
therefore relies in good part on the quality of the approximation of $u$ by
such simpler functions, which is typically governed
by the smoothness of $u$.

In our present setting, the fact that $u$
obeys the Helmholtz equation, may be used in addition
to its smoothness in order to guarantee the accuracy of certain approximation 
schemes, which are well adapted to such solutions. 

\begin{figure}
\centering
\includegraphics[width=4cm]{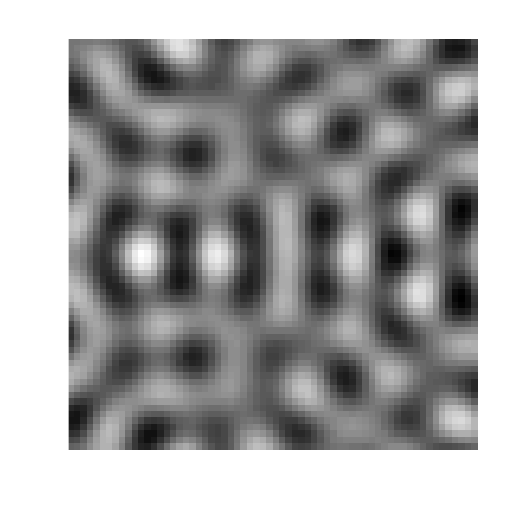}
\includegraphics[width=4cm]{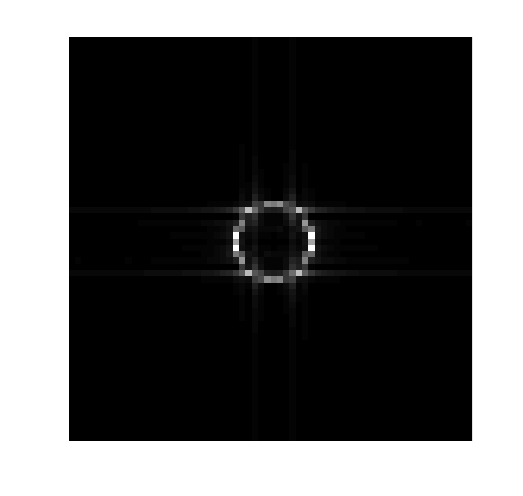}
\caption{A solution to the Helmholtz equation in 2D and its 2D Discrete Fourier Transform}
\label{modefft}
\end{figure}

The fact that the function to be measured is solution to the Helmholtz equation can
be used in various ways:
\begin{itemize}
\item As can be seen on Fig. \ref{modefft}, the spectrum of a solution to the equation (in 2D) is concentrated on an annulus.
This annulus can be enclosed in a square, or an hexagon, allowing reconstruction of the function from its values
on a square or hexagonal lattice, using the Shannon-Nyquist sampling theorem.
\item The recent field of compressed sensing suggests to interpret this property as the sparsity of the function
in a dictionary of Fourier-modes, and to reconstruct the function from a random sampling of the function on the
domain of interest.
\item The function can be reconstructed from its value on the border on the domain, as well as the value of its
normal derivative, using the Green formula:
$$ u(x) = \int_{\partial\Omega} u(y) \frac{\partial G}{\partial n}(y,x) - G(y,x) \frac{\partial u}{\partial n}(y) ds.$$
\end{itemize}

The last method is however not relevant for our setting, in which we are allowed to measure the function but not its derivatives.
The two first methods will be compared with the method we propose, based on the expansion of solutions to the Helmholtz equation
on particular families of function and least-squares approximations.

We introduce the Fourier-Bessel functions
\be
b_{\lambda, j}(x):= e^{ij\theta} J_j(\lambda r)
\label{fourierbessel}
\ee
where $(r, \theta)$ are the polar coordinates of $x$ and $J_j$ is the $j$-th Bessel function of the first kind. $b_{\lambda, j}$ is solution
to the Helmholtz equation \iref{helm} over $\R^2$ if and only if 
its parameter $\lambda$ is the same as in \eref{helm}.
Denoting $V_\lambda$ the set of the solutions of \eref{helm},
 it is known \cite{HMP2} that
 \be
V_\lambda=\overline{{\rm span}\{b_{\lambda, j}\}}^{L^2(\Omega)},
\ee
 and that the solutions of \eref{helm} can be approximated by elements of
the subspaces $V^b_{m} = {\rm span}\{b_{\lambda, j}, -L \leq j \leq L \}$
as $m = 2L+1$ grows.

An alternative approximation scheme uses plane waves defined by 
\be
e_{\bk}(x) := e^{i \bk \cdot x}
\ee
which are solutions of \eref{helm} if and only if $|\bk|=\lambda$.
The spaces $V^e_m$, spanned by the particular plane waves 
\be
e_{j}:=e_{\bk_j}, \;\; \bk_j:=\lambda\left(\cos\left(\frac {2j\pi}{m}\right),\sin\left(\frac {2j\pi}{m}\right)\right), \;\; j=-L,\dots,L,
\ee
can also be used to approximate solutions of \eref{helm} as $m$ grows \cite{HMP2}.

The most widely used approach to approximate $u$
in a finite dimensional space $V_m$, from its data at points
$x_1,\dots,x_n$, is the least squares method, namely
with $m\leq n$ solving the minimization problem
\be
\pi=\argmin_{v\in V_m}\frac 1 n\sum_{i=1}^n|y_i-v(x_i)|^2.
\label{ls}
\ee
The effectiveness of the least squares approximation
is governed by a certain trade-off in the choice of the
dimension $m$ of the approximation:
\begin{itemize}
\item 
A small value of $m$ leads to a highly regularized
reconstruction of $u$, which is usually robust but 
has poor accuracy.
\item
A large value of $m$ may lead to unstable
and therefore inaccurate reconstructions although the
space $V_m$ contains finer approximants to $u$.
\end{itemize}
Let us observe that regularization is relevant even
in a noiseless context where the function is measured exactly:
for example choosing $m=n$ corresponds to searching for
an exact interpolation of the data which may be very 
unstable and inaccurate, a phenomenon similar to the Runge phenomenon
in polynomial approximation.

In this paper, we discuss the amount of regularization
which is needed when applying the least squares method
using the finite dimensional subspaces $V^b_m$ and $V^e_m$
extracted from
$V_\lambda$.
With such discretizations, the distribution of the 
sampling points $x_1,\dots,x_n$ has an influence on the above described trade-off.
Our main theoretical result, established in the case
of a disc, shows that higher values of 
$m$, leading therefore to better accuracy, can be used if the $x_i$
are not uniformly distributed on $\Omega$ in the sense that a fixed fraction of these
points are located on the boundary $\partial \Omega$.
This result is confirmed by numerical experiments.

The rest of this paper is organized as follow: we give a brief account
in \S 2 on approximation of solutions to \iref{helm} by Fourier-Bessel functions and plane waves
which relies on Vekua's theory, and in \S 3 on general results
on the stability and accuracy of least-squares approximations
recently established in \cite{CDL}. We then study in \S 4 
the spaces $V^e_m$ and $V^b_m$ in more detail, in the
particular case where $\Omega$
is a disc, and use the above mentioned results to 
compare least-squares approximations on these spaces based on 
different sampling strategies. We also give similar results
for the case of the 3D ball.
In \S 5, we present numerical tests that illustrate the validity 
of this comparison. We also
show that our approach compares favorably with reconstructions
based on other approximation schemes such as partial Fourier sums 
(that do not exploit the fact that $u$ is a solution to \iref{helm}) 
and to other form of regularizations such as weighted Basis Pursuit.
In \S 6, we discuss further issues, namely determination of
the model order via cross-validation, and the influence of the
sampling distribution in the treatment of more general domains.

\section{Approximation by Fourier-Bessel functions and plane waves}

Results on the approximation of solutions to (\ref{helm}) 
by Fourier-Bessel functions and plane waves given in \cite{HMP2} are based on
the theory developed in the 1950's by Vekua \cite{Ve}. This theory generalizes
approximation results for holomorphic functions, 
viewed as solutions of $\Delta u=0$, to solutions of more general elliptic partial differential equations,
by means of appropriate operators that link the two types of solutions.

In the case of  the Helmholtz equation on a domain $\Omega$, that is star-shaped
with respect to a point which is fixed as the origin $0$, these operators (in their version
mapping harmonic functions to solutions to \iref{helm}) have
the explicit expression
\be
\cV_1\phi(x) = \phi(x) - \frac{\lambda |x|}{2} \int_0^1 \frac{1}{\sqrt{1-t}} J_1(\lambda |x| \sqrt{1-t}) \phi(tx)dt,
\ee
and
\be
\cV_2\phi(x) = \phi(x) - \frac{\lambda |x|}{2} \int_0^1 \frac{1}{\sqrt{t(1-t})} I_1(\lambda |x| \sqrt{1-t})\phi(tx) dt,
\ee
where $|x|$ stands for the euclidean norm of $x$, $J_1$ is the Bessel function of the first kind of order $1$ and $I_1$ the modified Bessel function of the first kind of order $1$, see \cite{HMP1} for more details.
These operators have important properties:
\begin{itemize}
\item They are linear.
\item $\cV_1$ maps harmonic functions to solutions of the Helmholtz equation,
and $\cV_2$ does the converse.
\item When restricted to harmonic functions or solutions to the Helmholtz equation, they are continuous in the Sobolev $H^k$ norms for all $k\geq 0$.
\item They are inverse to each other on these spaces.
\end{itemize}

As a consequence, any approximation method for harmonic functions can be translated as an approximation method for solutions of the Helmholtz equation.
In particular, approximation of harmonic functions by harmonic polynomials of degree $m$ translates as approximation of solution of the Helmholtz equation
by the so-called generalized harmonic polynomials $m$ which are
their image by $V_1$. The generalized
harmonic polynomials of degree $m$ can be expressed as
linear combinations of the Fourier-Bessel functions 
\iref{fourierbessel}, leading therefore to results
for the approximation of solutions to \iref{helm}
by elements of $V_m^b$ in Sobolev norms. More precisely, the following result can be obtained
when the domain $\Omega$ is convex, see theorem 3.2 of \cite{HMP2}
\be
\min_{v\in V^b_m} \|u - v\|_{H^k} \le C \left( \frac{\log m }{m}\right)^{p-k}
\|u\|_{H^p},
\label{approbes}
\ee
where the constant $C$ depends on $p$, $k$, $\lambda$ and the geometry of $\Omega$.
This results still holds for more general, star-shaped convex domains,
with a slower convergence.

Plane waves and Bessel functions are related by the Jacobi-Anger identity
\be
e_{\mathbf k_\phi} = \sum_{m\in\Z} i^m J_m(\lambda r) e^{im(\theta - \phi)}.
\ee
where $\bk_\phi:=\lambda(\cos(\phi),\sin(\phi))$, and its converse,
the Bessel integral
\be
J_n(\lambda r) e^{in\theta} = \frac{1}{2 \pi i^n} \int_{-\pi}^{\pi} e_{\bk_\phi} e^{in\phi} d\phi,
\label{jacobiangerconv}
\ee
Approximating the integral in \iref{jacobiangerconv} by a discrete sum,
by uniformly sampling the wave vectors $\bk_\phi$ on the circle of diameter $\lambda$, 
leads to approximations of solutions to \iref{helm} by linear combinations of the $2m+1$ plane waves $e_{\mathbf k_j}$,
that is, by elements of $V_m^e$. It is also known (see theorem 5.2 of \cite{HMP2})
 that
such approximations have the same convergence properties, e.g., for a convex
domain,
\be
\min_{v\in V^b_e} \|u - v\|_{H^k} \le C \left( \frac{\log m }{m}\right)^{p-k}
\|u\|_{H^p}.
\label{approplan}
\ee

\section{Least-square approximations}

The results of the previous section quantify how 
a general solution $u$ to the Helmholtz equation can be approximated
by functions from spaces $V^e_m$ or $V^b_m$. We are now interested in 
understanding the quality of approximations
from these spaces built by the least squares methods
based on scattered data $(x_l,y_l)_{l=1,\dots,n}$.
In particular, we want to understand the 
trade-off between the dimension $m$
and the number of samples $n$. Ideally we 
would like to choose $m$ large in order to benefit
of the approximation properties \iref{approbes} and \iref{approplan}, however not too large so that
stability of the least-square method is ensured.
We are also interested in understanding how the spatial
distribution of the sample $x_l$ influences this trade-off.

This problem was recently studied in \cite{CDL}, 
in a general setting where the $x_l$ are independently drawn
according to a given probability measure
$\nu$ defined on $\Omega$. This measure therefore
reflects the spatial distribution of the samples. For example,
the uniform measure 
\be
d\nu:=|\Omega|^{-1} dx,
\ee
tends to generate uniformly spaced samples. In order to present
the general result of \cite{CDL},
we assume that $(V_m)_{m\geq 1}$ is an arbitrary sequence 
of finite dimensional spaces of functions defined on $\Omega$ with $\dim(V_m)=m$.

We introduce the $L^2$ norm with respect to the measure $\nu$
\be
\|v\|:=\left(\int_\Omega |v|^2 d\nu\right)^{1/2},
\ee
and we define the best approximation error for a function $u$ in this norm as
\be
\sigma_m(u):=\min_{v\in V_m} \|u-v\|.
\ee
Note that, in the noiseless case, the least squares method amounts to 
computing the best approximation of $u$ onto $V_m$ with respect to the 
norm
\be
\|v\|_n:=\left(\frac 1 n\sum |v(x_l)|^2\right)^{1/2}.
\ee
This norm can be viewed as an approximation of the norm $\|v\|$ based on the draw,
and it is therefore natural to compare the
error $\|u-\pi\|$ where $\pi$ is computed by \iref{ls}
with $\sigma_m(u)$. We give below a 
criterion that describes under which condition
on $m$ these two quantities are of comparable size.

 Here, we assume that $(L_1,\dots,L_m)$ is a basis of $V_m$
 which is orthonormal in $L^2(\Omega,\nu)$. We define the quantity
 \be
 K(m)=K(V_m,\nu):=\max_{x\in\Omega}\sum_{j=1}^m |L_j(x)|^2,
 \ee
 which depends both on $V_m$ and on the chosen measure $\nu$,
 but not on the choice of the orthonormal basis since it is invariant by rotation.

 We also assume that an a-priori bound $\|u\|_{L^\infty}\leq M$ is known
 on the function $u$. We can therefore only improve the
 least squares estimate by defining
 \be
 \tilde u:=T_M(\pi),
 \ee
 where $T_M(t):={\rm sign}(t)\min\{|t|,M\}$ and $\pi$ is given by \iref{ls}.
 The following result was established in \cite{CDL}, in the case of noiseless 
 data, i.e.\ $\eta_l=0$ in \iref{data}.
 
 \begin{theorem}
 Let $r>0$ be arbitrary but fixed and let $\kappa:=\frac {1-\log 2}{2+2r}$. 
 If $m$ is such that 
 \be
 \label{condkm}
 K(m)\leq \kappa \frac n {\log n},
 \ee
 then, the expectation of the reconstruction error is bounded:
 \be
 E(\|u-\tilde u\|^2)\leq (1+\e(n))\sigma_m(u)^2+8M^2n^{-r},
 \ee
 where $\e(n):=\frac {4\kappa} {\log n} \to 0$ as $n\to +\infty$.
\label{theols}
\end{theorem}

It is also established in \cite{CDL} that the condition \iref{condkm}
ensures the numerical stability of the least-square method, with probability
larger than $1-2n^{-r}$.
These results suggest to set the regularization level
by picking the largest value of $m^*=m^*(n)$
such that \iref{condkm} holds. The dependence of $m^*(n)$ with $n$ is 
obviously related to that of $K(m)$ with $m$. In particular,
slower growth of $K(m)$ with $m$ implies faster growth
of $m^*(n)$ and therefore faster convergence of the least squares approximation.
Notice that we always have 
\be
K(m) \geq \int_\Omega \sum_{j=1}^m |L_j|^2 d\nu=m.
\ee
In the next section, we evaluate $K(V^e_m,\nu)$ and $K(V^e_b,\nu)$ in the case where
$\Omega$ is a disk, for various choices of the measure $\nu$.

\section{Stability of the reconstruction on a disc and in a ball}

\begin{figure}
\centering
\includegraphics[width=8cm]{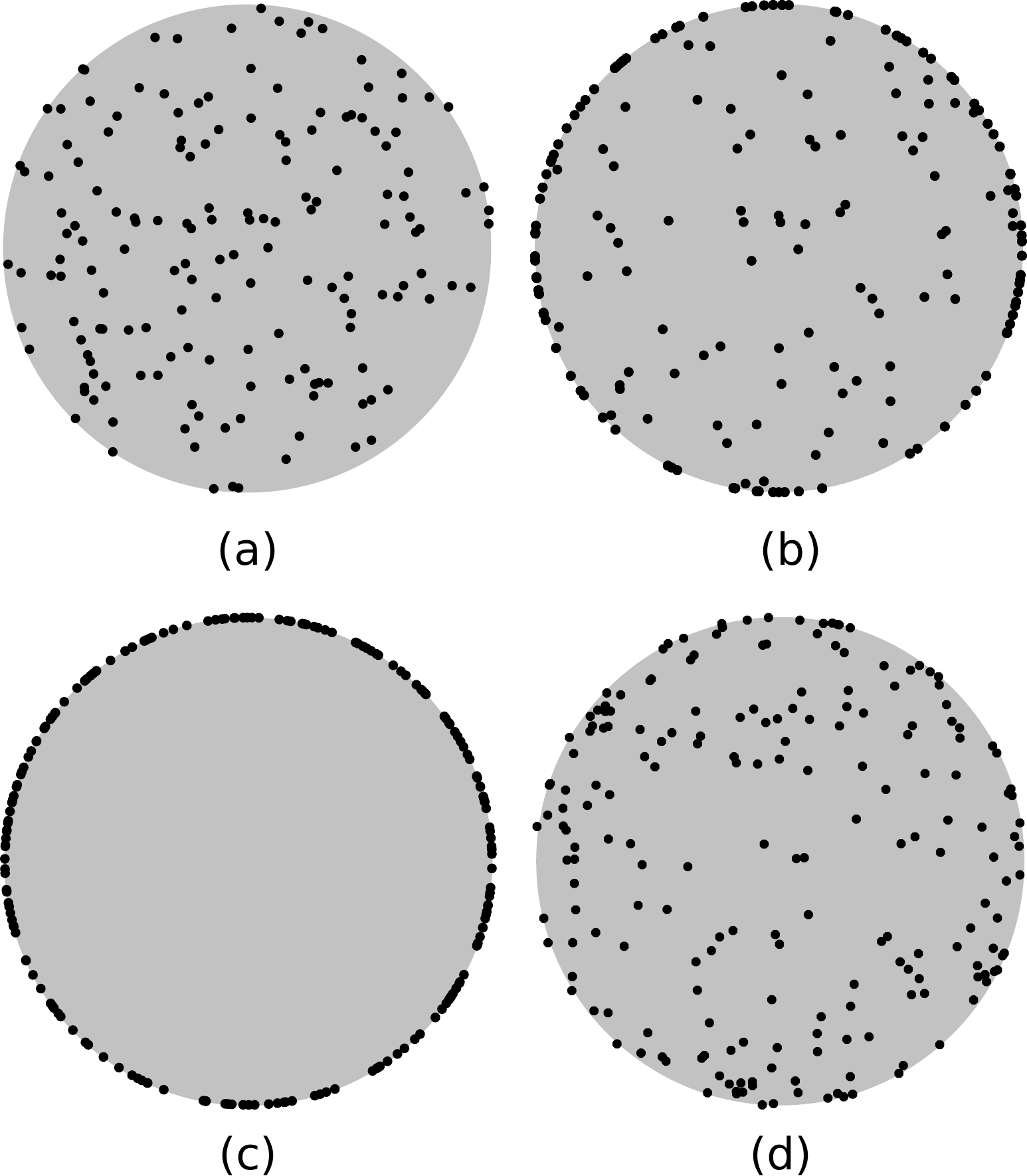}
\caption{Examples of sampling distributions on the disc:
(a) $\nu_0$, (b) $\nu_{1/2}$, (c) $\nu_1$ (d) $\nu'$, defined and used in section \ref{num}.}
\label{circlespl}
\end{figure}

As mentioned in the previous section, the quantity $K(m)$ depends
both on the space $V_m$ and the measure $\nu$ that reflects
the sampling strategy. Here we study the case where $V_m$ is either
one of the spaces of plane waves $V^e_m$ or of Fourier-Bessel functions $V^b_m$ defined in the introduction,
for $m = 2L+1$.
We consider two sampling strategies. The first one uses the uniform probability distribution
\be
\nu_0:=\frac {dx}{|\Omega|},
\ee
and the second one combines uniform sampling on the domain and on its boundary, with proportion 
$0<\alpha < 1$, according to the probability distribution
\be
\nu_\alpha :=(1-\alpha)\frac {dx}{|\Omega|}+\alpha \frac {d\sigma}{|\partial\Omega|}.
\label{nualpha}
\ee
Examples of such densities are pictured on figure \ref{circlespl}. The norm computed using these densities is simply denoted $\|\cdot\|$.

This particular choice of probability distribution makes it possible to control
the $L_2(\Omega)$-norm of the error. Indeed, theorem \ref{theols} control the reconstruction
error in the norm defined by $\nu_\alpha$, which itself can bound the $L_2(\Omega)$-norm of the error as
\be
\|u\|_{L_2(\Omega)} \leq \frac{1}{1-\alpha} \|u\|.
\label{bound}
\ee

In order to obtain explicit results, we focus on the simple
case where $\Omega$ is a disk. Without loss of generality, we fix
\be
\Omega:=\{x\in\R^2 \; : \; |x|\leq 1\}.
\ee

\subsection{Fourier-Bessel approximation on the disc}

Fourier-Bessel functions are orthogonal on the disk with respect to any rotation-invariant measure, because of their angular
dependence in $e^{in\theta}$. This allows a simple computation of the quantity $K(m)$ for the space $V_m^b$,
leading to the following result.

\begin{theorem}
\label{Theofourbes}
For the space $V_m^b$ and the measure $\nu_\alpha$ on the unit disk $\Omega$, one has
for sufficiently large $m$
\be
K(m)\geq c_0 + c_1 m^2,
\label{res11}
\ee
when $\alpha=0$ (that is, for the uniform measure), for any $c_1 < 1/16$, and where $c_0$ depends on $c_1$ and $\lambda$, and
\be
K(m) \leq C+\frac {m} \alpha,
\label{res12}
\ee
when $\alpha>0$, where $C$ depends on $\lambda$ and $\alpha$.
\end{theorem}

This result indicates that using an order $m$ for the approximation necessitates a number of samples $n$ that scales at least quadratically
with $m$ when sampling uniformly in the disk. Using a proportion $\alpha$ of samples
on the border makes $K(m)$ linear with respect to $n$. The best behavior possible for $K$ (i.e.\ $K(m) = m$) can be approached when $\alpha$ approaches
1. However in that case, the constant $C$ may grow,
and the bound \eref{bound} becomes less efficient. This would make the use of a large proportion of samples on the border relevant
only for very large numbers of samples.

Note finally that in the case $\alpha=1$, Eq. \eref{bound} cannot be used to control the $L_2(\Omega)$-norm of the error,
allowing arbitrary large errors with any number of samples. For instance, when
$\lambda$ is an eigenfrequency of the disk with Dirichlet boundary conditions, the associated eigenmode
(a Fourier-Bessel function) cannot be recovered as its samples on the border are identically zero.

\noindent
{\bf Proof:}
Since the Fourier-Bessel functions are orthogonal in $L^2(\Omega,\nu_\alpha)$,
we have
\be
K(m) = \left \| \sum_{j=-L}^{L} \frac{|b_j|^2}{\|b_j\|^2}\right \|_{L^\infty(\Omega)}.
\ee
In the case of the uniform measure, we bound $K(m)$ from below.
We first write
\be
K(m) \geq  \left \|  \sum_{j=-L}^{L} \frac{|b_j|^2}{\|b_j\|^2}  \right\|_{L^\infty(\partial\Omega)}.
 \ee
We next bound $\|b_j\|^2$, for $j> \lceil \lambda \rceil$ (the case $j < -\lceil \lambda \rceil$
is identical, as $|b_j| = |b_{-j}|$), according to
\begin{eqnarray*}
\|b_j\|^2 & = & \frac{1}{\pi} \int_\Omega J_j(\lambda |x|)^2 dx \\
& = & 2 \int_0^1 r J_{j}(\lambda r)^2 dr \\
& = & \frac{4}{\lambda^2} \sum_{p=0}^\infty (j + 1 + 2p) J^2_{j+1+2p}(\lambda)\\
& \leq  & \frac{4}{\lambda^2} \sum_{p=0}^\infty (j + 1 + 2p) (\lambda/j)^{2+4p}J^2_{j}(\lambda)\\
& = & 4 \frac{j+1}{j^2}\left(\frac{1}{1 -(\lambda/j)^4} + \frac{2}{j+1}\frac{(\lambda/j)^4}{\left(1 - (\lambda/j)^4\right)^2} \right)J^2_{j} (\lambda) 
\end{eqnarray*}
where the third equality is identity (11.3.32) of \cite{AS}, and the first inequality comes from (A.6) of \cite{PD}.
As $|b_j(1,\theta)| = |J_{j}(\lambda)|$, we have, for any $c <1$
and $j$ larger than some  $j_0$
\be \frac{|b_j(1,\theta)^2|}{\|b_j\|^2} \ge 
\frac{j^2}{4(j+1)} \left(\frac{1}{1 -(\lambda/j)^4} + \frac{2}{j+1}\frac{(\lambda/j)^4}{\left(1 - (\lambda/j)^4\right)^2} \right)^{-1}
\ge c \frac{j}{4}
\ee
and
\begin{eqnarray}
K(m) & \geq & \sum_{j=-j_0}^{j_0} \frac{|J_j(\lambda)|^2}{\|b_j\|^2} + 2c \sum_{j_0 < j \leq L} \frac{j}{4}\\
& \geq & \sum_{j=-j_0}^{j_0} \frac{|J_j(\lambda)|^2}{\|b_j\|^2} - c \frac{j_0(j_0+1)}{4}
+ \frac{c}{4} L(L+1)
\end{eqnarray}
which proves the bound \iref{res11} when $L > j_0$.

In the case of mixed sampling, we bound $K(m)$ from above  by
\be
K(m)  \leq  \sum_{j=-L}^{L}  \frac{\|b_j\|_{L^\infty(\Omega)}^2}{\|b_j\|^2}.
\ee

$\|b_j\|^2$ is nonzero as $\|b_j\|^2 > (1-\alpha) \|b_j\|_{L_2(\Omega)}^2$ and $\alpha <1$.
When $j > \lambda$, the function $r\mapsto J_l(\lambda r)$ is monotone increasing on $[0,1]$, so that
$\|b_j\|_{L^\infty(\Omega)} = J_{j}(\lambda)$. Thus,
\begin{eqnarray*}
\frac{\|b_j\|_{L^\infty(\Omega)}}{\|b_n\|^2} & = & \frac{J_j(\lambda)^2} { \frac{1 - \alpha}{\pi} \int_D J_j(\lambda |x|)^2 dx
+ \frac{\alpha}{2\pi} \int_0^{2\pi} J_j(\lambda)^2 d\theta} \\
& < & \frac{1}{\alpha}.
\end{eqnarray*}
We then have $$K(m) \leq \frac{2(L - \lfloor \lambda \rfloor)}{\alpha} +  \sum_{j=-\lfloor \lambda \rfloor}^{\lfloor \lambda \rfloor} \frac{\|b_j(x)\|^2_{L^\infty(\Omega)}}{\|b_j\|^2}$$
which proves \iref{res12}.

\subsection{Plane wave approximation on the disc}

Similar results can be obtained for the plane wave approximation:
\begin{theorem}
For the space $V_m^e$ and the measure $\nu_\alpha$ on the unit disk $\Omega$, one has
for sufficiently large $m$
\be
K(m)\geq c_0 + c_1m^2,
\label{res21}
\ee
when $\alpha=0$ (that is, for the uniform measure), for any $c_1 < 1/16$, where $c_0$ depends on $c_1$ and $\lambda$, and 
\be
K(m) \leq C_1+ C_2\frac {m} \alpha,
\label{res22}
\ee
when $\alpha>0$, for any $C_2 > 1$, where $C_1$ depends on $C_2$, $\lambda$ and $\alpha$.
\end{theorem}

\noindent
{\bf Proof:} Since plane waves are not orthogonal in $L^2(\Omega,\nu)$,
the first step is the computation of an orthogonal basis
of the space spanned by these plane waves. Let us consider $2L+1$
plane waves, with wave vectors uniformly distributed on the circle
of radius $\lambda$. For the measures considered here, the Gram matrix of this family is a circulant matrix,
which is diagonalized in the Fourier basis. An orthogonal family
spanning the same space is therefore given by functions that are linear combinations
of plane waves with the coefficients of the discrete Fourier transform:
\be
b^m_j := \frac{1}{m} \sum_{j=-L}^L e^{2\pi i j/m} e_{\bk_j}.
\label{sumwaves}
\ee
This formula may be thought as a quadrature for the Bessel integral 
\iref{jacobiangerconv}: the
$b^m_j$ are thus approximations of the Fourier-Bessel functions $b_j$.
In order to bound the quantity $K(m)=K(V_m^e,\nu)$, we compare
it to the quantity $K(V^b_m,\nu)$ which behavior is described by Theorem \ref{Theofourbes}.
Using (8) and (9) from \cite{PD} we have
\be
b^m_j = \sum_{p\in \Z} i^{pm} b_{j+pm},
\label{aliasbessel}
\ee
and
\be
\left|\frac{b^m_j}{b_j} -1 \right| \leq \sum_{p\in \Z-\{0\}} \frac{|b_{j+pm}|}{|b_j|}.
\ee
With $L \geq j \geq \lambda$, we thus have for all $0\leq r\leq 1$ and $0\leq \theta\leq 2\pi$,
\begin{eqnarray*}
\left|\frac{b^m_j(r,\theta)}{b_j(r,\theta)} -1 \right| & \leq & \frac{1}{|b_j(r,\theta)|} \left(\sum_{p\geq 0} |b_{j+(p+1)m}(r,\theta)| + \sum_{p\geq 0} |b_{j+(p+1)m- 2j}(r,\theta)| \right)\\ 
& \leq & \sum_{p\in \N} \left(1 + \left(\frac{\lambda}{j}\right)^{-2j}\right) \left(\frac{\lambda}{j}\right)^{(p+1)m} \\
&= & \frac{\left(\frac{\lambda}{j}\right)^m + \left(\frac{\lambda}{j}\right)^{m-2j}}{1 - \left(\frac{\lambda}{j}\right)^m} \\
& \leq &  \frac{2\frac{\lambda}{j}}{1 - \frac{\lambda}{j}}
\end{eqnarray*}
where we have used equation (A.6) of \cite{PD} to obtain the second inequality.

Using the orthogonality of the $b_j$, we have
\be
\|b^m_j\|^2= \sum_{p\in \mathbf Z} \|b_{j+pm}\|^2,
\ee
and
\be
\left| \frac{\|b^m_j\|^2}{\|b_j\|^2} -1 \right| = \sum_{p\in \Z-\{0\}} \frac{\|b_{j+pm}\|^2}
{\|b_j\|^2}.
\ee
When $j \geq \lambda$ and $l \geq 0$ we bound $\|b_{j+l}\|^2$ according to
\begin{eqnarray*}
\|b_{j+l}\|^2  &= & 2\pi \int_{0}^1 r J_{j+l}(\lambda r)^2 dr \\
& \leq & 2\pi \int_{0}^1 r \left(\frac{\lambda r}{j} \right)^{2l} J_j(\lambda r)^2 dr \\
& \leq &\left(\frac{\lambda}{j} \right)^{2l}2\pi \int_{0}^1 J_j(\lambda r)^2 dr\\
& = &\left(\frac{\lambda}{j} \right)^{2l} \|b_{j}\|^2
\end{eqnarray*}
where we again have used equation (A.6) of \cite{PD} to obtain the first inequality.
We thus have
\begin{eqnarray*}
\left| \frac{\|b^m_j\|^2}{\|b_j\|^2} -1 \right| & = & \sum_{p\in \Z-\{0\}} \frac{\|b_{j+pm}\|^2}{\|b_j\|^2} \\
& = & \frac{1}{\|b_j\|^2} \left(\sum_{p\geq 0}  \|b_{j+(p+1)m}\|^2  +  \sum_{p\geq 0}\|b_{j+(p+1)m -2j}\|^2\right) \\
& \leq & \left(\sum_{p\geq 0}  \left(\frac{\lambda}{j}\right)^{2(p+1)m}  + \sum_{p\geq 0} \left(\frac{\lambda}{j}\right)^{2((p+1)m -2j)}\right) \\
& = & \frac{\left(\frac{\lambda}{j}\right)^{2m} + \left(\frac{\lambda}{j}\right)^{2(m-2j)}}{1 - \left(\frac{\lambda}{j}\right)^{2m}} \\
& \leq & \frac{2\left(\frac{\lambda}{j}\right)^{2} }{1 - \left(\frac{\lambda}{j}\right)^{2}}. \\
\end{eqnarray*}

We now consider the case $\alpha=0$. Using the above comparison results
between $b_j$ and $b_j^m$, we can find, for any $c<1$ a positive integer $j_c$ such that for $L \geq j \geq j_c$,
\be  
\frac{|b^m_j(1,\theta)|^2}{\|b^m_j\|^2} \geq c  \frac{J_j(\lambda)^2}{\|b_j\|^2}.
\ee
We may thus write
\begin{eqnarray*}
K(V_m^e,\nu) & = & \left \| \sum_{j=-L}^{L} \frac{|b_j^m|^2}{\|b_j^m\|^2} \right \|_{L^\infty(\Omega)} \\
& \geq & \sum_{j=-L}^{ L} \frac{|b_j^m(1,0)|^2}{\|b_j^m\|^2} \\
& \geq & \sum_{j < j_c} \frac{|b_j^m(1,0)|^2}{\|b_j^m\|^2} + 2 c  \sum_{j=j_c}^{L} \frac{J_j(\lambda)^2}{\|b_j\|^2},\\
& \geq & 2 c  \sum_{j=j_c}^{L} \frac{J_j(\lambda)^2}{\|b_j\|^2}.\\
\end{eqnarray*}
The last sum can be bounded from below in a similar way as in the proof of Theorem \ref{Theofourbes}, proving \eref{res21}.

We next consider the case $\alpha>0$. We then write
\be
K(m)  \leq  \sum_{j=-L}^{L}  \frac{\|b_j^m\|^2_{L^\infty(\Omega)}}{\|b_j^m\|^2}.
\ee
For any $C>1$, there is a $j_C > \lambda$ such that when $j > j_C$, we have $\frac{1}{C} |b_j(r,\theta)| \leq |b_j^m(r,\theta)| \leq C |b_j(r,\theta)|$, so that
$ \|b_j^m\|_{L^\infty(\Omega)} < C \|b_j\|_{L^\infty(\Omega)} = C J_j(\lambda)$, and $|b^m_j(1,\theta)| \geq  |J_j(\lambda)|/C$.
We then have, for $m \geq j \geq j_C$,
\begin{eqnarray*}
\frac{\|b_j^m\|^2_{L^\infty(\Omega)}}{\|b^m_j\|^2} & = & \frac{\|b_j^m\|^2_{L^\infty(\Omega)}}{ \frac{1 - \alpha}{\pi} \int_D |b_j^m|^2 dx
+ \frac{\alpha}{2\pi} \int_0^{2\pi} |b_j^m|^2 d\theta}\\
& \leq & \frac{1}{\alpha} \frac{C^2 J_j(\lambda)^2}{J_j(\lambda)^2 /C^2}\\
& \leq &\frac{C_2}{\alpha}
\end{eqnarray*}
and
\begin{eqnarray*}
K(m) & \leq & \frac{2C_2(L - \lfloor \lambda \rfloor)}{\alpha} +  \sum_{|j| < j_0}  \frac{\|b_j^m\|^2_{L^\infty(\Omega)}}{\|b_j^m\|^2} \\
& \leq & \frac{2C_2(L - \lfloor \lambda \rfloor)}{\alpha} + \sum_{|j| < j_0}  \frac{1}{\|b_j\|^2}
\end{eqnarray*}
which proves \ref{res22}. We use here the fact that $\|b^m_j\|_{L^\infty(\Omega)} \leq 1$ which is clear from
\eref{sumwaves}, and $\|b_j\| \leq \|b_j^m\|$ obtained from \eref{aliasbessel} and the orthogonality
of the $b_j$. \hfill $\Box$

\subsection{Spherical Fourier-Bessel functions in a ball}

Similar results can be obtained for the approximation in a ball.
 In the 3D case, solutions to the Helmholtz equation can be approximated by sums of products of spherical harmonics $Y_{l,q}$ and spherical Bessel functions $j_l$ \cite{HMP1}:
$$b_{\lambda, l,q}(x) = Y_{l,q} \left(\frac{x}{|x|}\right) j_l (\lambda |x|).$$
For $m= (L+1)^2$, the $(L+1)^2$-dimensional space $V_m^3$ is defined as $V^3_{m} = {\rm span}\{b_{\lambda, l, q}, 0 \leq q \leq L, -l \leq q \leq l \}$,
and
\be
\min_{v\in V^3_m} \|u - v\|_{H^k} \le C m^{-\alpha(p-k)}
\|u\|_{H^p},
\label{approbes3}
\ee
where $\alpha$ is a strictly positive constant (in general, this constant depends on the
shape of the domain of interest).

Using similar sampling densities $\nu_\alpha$ (that is, a proportion $\alpha$ of the samples on the sphere, the rest inside the ball), we have:

\begin{theorem}
\label{Theofourbes3}
For the space $V_m^3$ and the measure $\nu_\alpha$ on the unit ball $\Omega$, one has
for sufficiently large $m = (L+1)^2$
\be
K(m)\geq c_0 + c_1 m^{3/2},
\ee
when $\alpha=0$ (that is, for the uniform measure), for any $c_1 < 1/9$, and where $c_0$ depends on $c_1$ and $\lambda$, and
\be
K(m) \leq C+\frac {m} \alpha,
\ee
when $\alpha>0$, where $C$ depends on $\lambda$ and $\alpha$.
\end{theorem}

{\bf Proof:} The proof is a straightforward adaptation of the proof of theorem \ref{Theofourbes}, using
$\sum_{q = -l}^l \left|Y_{l,q}\left(\frac{x}{|x|}\right)\right|^2 = 2l+1$
and
$j_l (t)= \sqrt{\pi/(2t)} J_{l+1/2}(t)$.

\section{Numerical tests}
\label{num}

Here, we compare four different reconstruction
methods on the unit disc:
\begin{enumerate}[(i)]
\item Least-squares method with a dictionary of Fourier modes.
\item Weighted $\ell_1$-minimization \cite{rauhut} with a Fourier dictionary.
\item The proposed method, least-squares method with a dictionary of Fourier-Bessel functions.
\item Weighted $\ell_1$-minimization with a Fourier-Bessel dictionary.
\end{enumerate}

For the first two methods, the dictionary contains orthogonal
Fourier modes on a square enclosing the disc. These modes thus have the form
$e^{i a \bk \cdot x}$ for some fixed $0<a\leq \pi$ (here we took $a=\pi/2$) and $\bk \in \{-K,\dots,K\}^2$. The size of the dictionary is smaller than the number of measurements for methods (i) and (iii), and larger for methods (ii) and (iv).

The methods are tested for $\lambda = 12$ with solutions that are the
linear combinations of fundamental solutions (i.e.\ second kind Bessel function
$Y_0(\lambda r)$ where $r$ is the distance to the source). The sources are placed
on a circle of radius 1.1. This setup can occur when synthesizing acoustical fields.

\begin{figure}
\centering
\includegraphics[width=10cm]{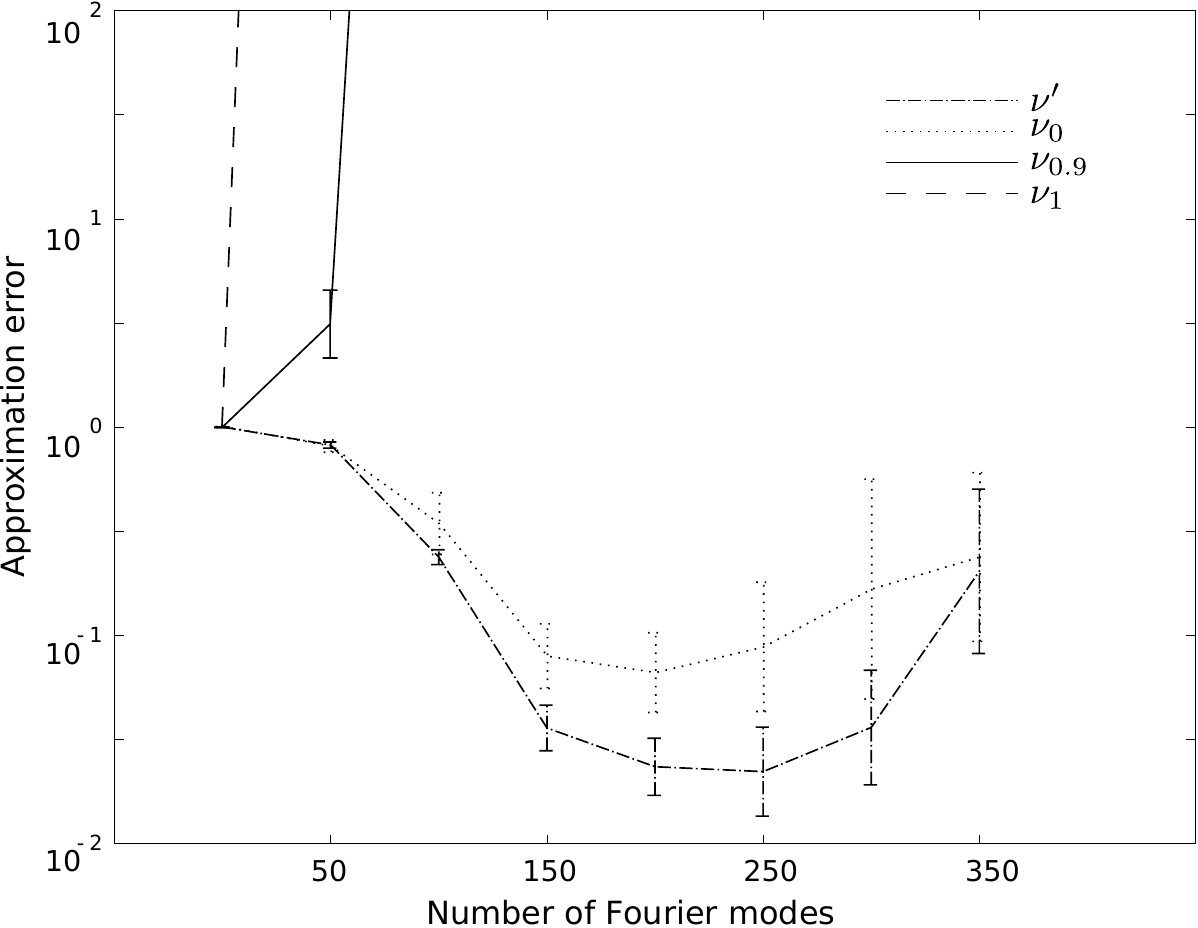}
\caption{Reconstruction error for method (i), least-squares with Fourier dictionary, vs. number of Fourier modes with $n=400$ measurements}
\label{ef400}
\end{figure}

\begin{figure}
\centering
\includegraphics[width=10cm]{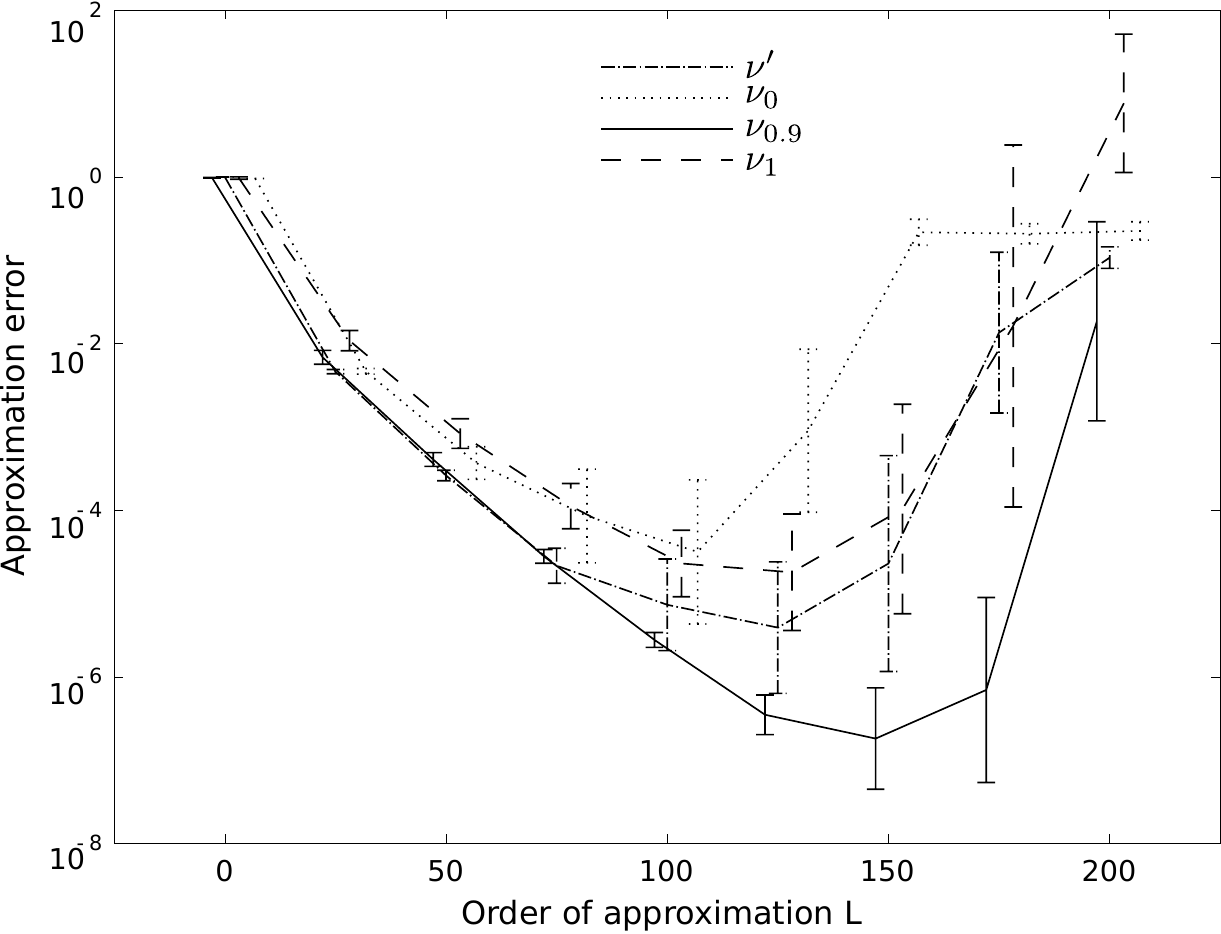}
\caption{Reconstruction error for the proposed method (iii), least-squares with Fourier-Bessel dictionary, vs. number of Fourier-Bessel functions with $n=400$ measurements}
\label{efb400}
\end{figure}

Results for the least-square methods (i) and (iii) are given on figures
\ref{ef400} and \ref{efb400}, for the different values $\alpha 
=0, 0.9, 1$. Another distribution $\nu'$ is also tested,
with uniform distribution in angle, and radiuses drawn from the interval $(0, 1)$ with probability density $\pi/ 2\sqrt{1 - r^2}$. An example of such
distribution is given on Fig. \ref{circlespl}, showing the higher
density of samples near the boundary. We plot the error measured in the norm $L^2(\Omega)=L^2(\Omega,dx)$, averaged over 40 realizations of the sampling, versus the approximation space dimension $m$.

Reconstruction errors for method (i), with a Fourier dictionary, are
always above $10^{-2}$ and do not benefit from sampling on the boundary,
as the best results are obtained with $\nu_0$ (uniform sampling) or $\nu'$.

Results for the proposed method (iii), with the Fourier-Bessel dictionary, are displayed on Figure \ref{efb400}.
We observe that placing more measurement points on the boundary $\partial \Omega$ is beneficial
to the reconstruction:  for $\alpha = 0.9$, the error is reduced by four order of magnitude compared to method (i).

However, measuring
the solutions only on $\partial \Omega$ ($\alpha = 1$) does not yield good reconstructions.
In that case, Theorem \ref{theols}  only ensures that the reconstruction
 is accurate on $\partial \Omega$ and says nothing
on the error on the disk itself, since we cannot control the $L^2(\Omega)$ norm
by the $L^2(\Omega,\nu_1)$-norm, which is actually the $L^2(\partial\Omega)$-norm.

\begin{figure}
\centering
\includegraphics[width=10cm]{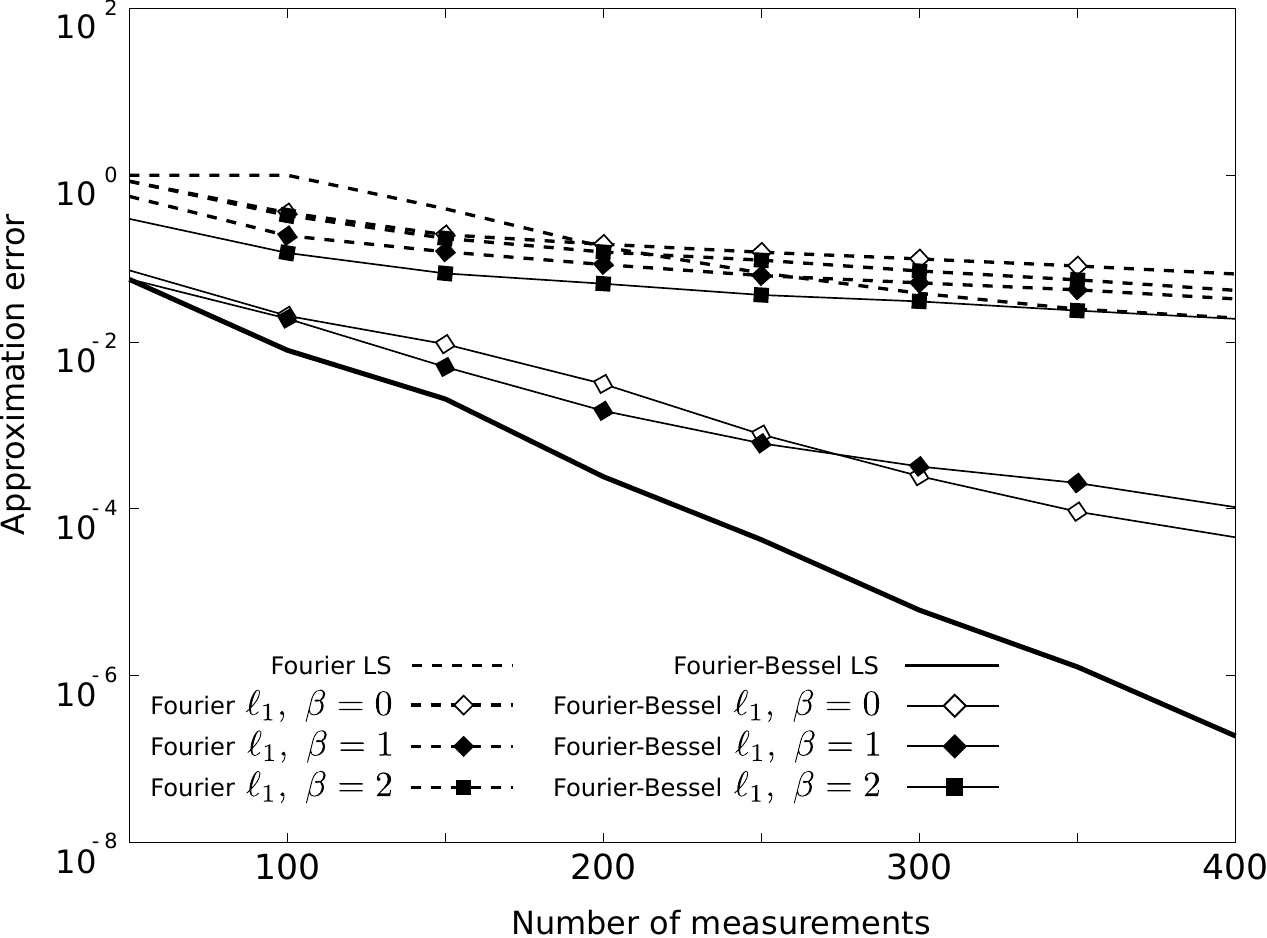}
\caption{Best reconstruction error for method (i) Fourier least-squares, (ii) Fourier Weighted Basis Pursuit, (iii) Fourier-Bessel least-squares (proposed method) and (iv) Fourier-Bessel Weighted Basis Pursuit}
\label{ecomp}
\end{figure}

Figure \ref{ecomp} compares the behavior of these methods ``at their best'' with varying number of
measurements, from 50 to 400. The plotted errors are obtained by selecting
the value of $K$ for (i), and of $m$ for (iii)
as well as the proportion $\alpha$, that minimize the error for the given number $n$ of measurements. Results 
for (i) and worse than for the proposed
method (iii), and are always obtained with the distribution $\nu'$ for the least-squares method with Fourier modes.
In contrast, as expected, the best results of (iii) are obtained with $\alpha = 0.9$.

Results of method (ii), Basis Pursuit with Fourier dictionary, are also given, with weights $(1+k)^{\beta}$ where
$k$ is the wavenumber of the Fourier mode, taking in consideration the
sparsity as well as the smoothness of the functions to be reconstructed. We use here the
SPGL1 toolbox \cite{spgl1, spgl1paper}.
Best results
are obtained for the sampling with the distribution $\nu'$ and
$\beta = 1$.  Performances of this
method are not as good as the Fourier-Bessel least-squares method. Using method (iv), i.e. the same algorithm with a larger Fourier-Bessel dictionary
that the one used for method (iii) yields, for $\beta = 0$ and $\beta = 1$
(the weight is here $(1+j)^\beta$ where $j$ is the order of the Fourier-Bessel function),
good reconstructions, but not as good as the simpler least-squares method.
Best results are here obtained for $\alpha = 0.9$.

Another sparse approximation method, Orthogonal Matching Pursuit \cite{omp} was also tested.
The reconstruction errors were always larger than the results of Basis Pursuit.

The better results obtained with the proposed method (iii) show that if an adequate
model is used to describe the signals of interest (here Fourier-Bessel 
approximations, capturing the particular type of sparsity exhibited by
the solutions to the Helmholtz equation better than a simpler Fourier dictionary) and an appropriate sampling scheme is used, basic numerical methods (here, standard least-squares estimation) yield better results than more sophisticated
methods such as weighted $\ell_1$-minimization.

\section{Further issues}

We showed above that a careful choice of the sampling distribution allows
us to use a larger order of approximation. However, in practice, the optimal value of $m$ is unknown to us.
We test here the  cross-validation  method to estimate this value.
We also estimate the value of $K(m)$ for a square and different sampling densities, and show that
a non-uniform density on the border may be necessary in a general setting.

\subsection{Estimation of the model order with the  cross-validation method}

\begin{figure}
\centering
\includegraphics[width=10cm]{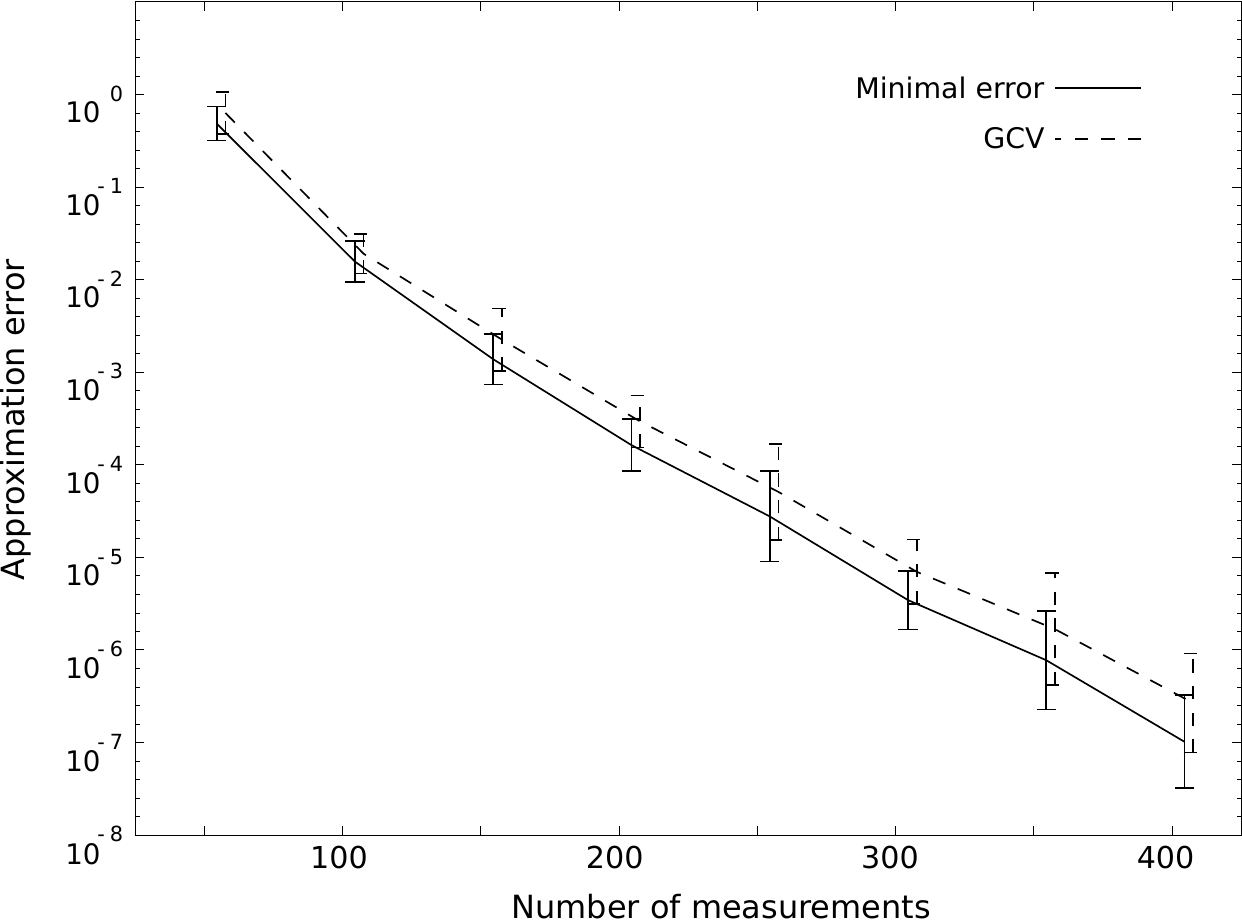}
\caption{Best reconstruction error for the proposed method, using GCV and optimal value of $m$.}
\label{egcv}
\end{figure}

For a given number $m = 2L+1$ of Fourier-Bessel functions,
we estimate a reconstruction using 95\% of the samples, and evaluate empirically the
mean-square error using the remaining ones. We repeat the estimation
10 times using different choices of estimation and reconstruction
points within the same sample, and select the number 
$m$ that minimizes the mean
square error. The function is then reconstructed using all samples. On
Figure \ref{egcv}, we compare the results obtained by this method 
with those based on the optimal value of $m$, as 
the number of measurement $n$ varies.
Here we use the sampling distribution according to the measure
$\nu_\alpha$, with $\alpha = 0.9$.
We observe that the performances are comparable
up to a slight loss by a multiplicative constant.

\subsection{More general shapes}

While the results obtained here inform us on the importance of sampling on the border of the considered domain,
a numerical test on another simple shape shows that the density on the border is critical.
The theoretical analysis for the disk and the ball was based on the fact that the Fourier-Bessel functions
were already an orthogonal basis of the space $V_m^b$. We focus here on the square $[-1, 1]^2$. As neither the Fourier-Bessel
functions, nor the plane waves, form an orthogonal basis, we construct
one by orthogonalizing the Fourier-Bessel functions.

\begin{figure}
\centering
\includegraphics[width=10cm]{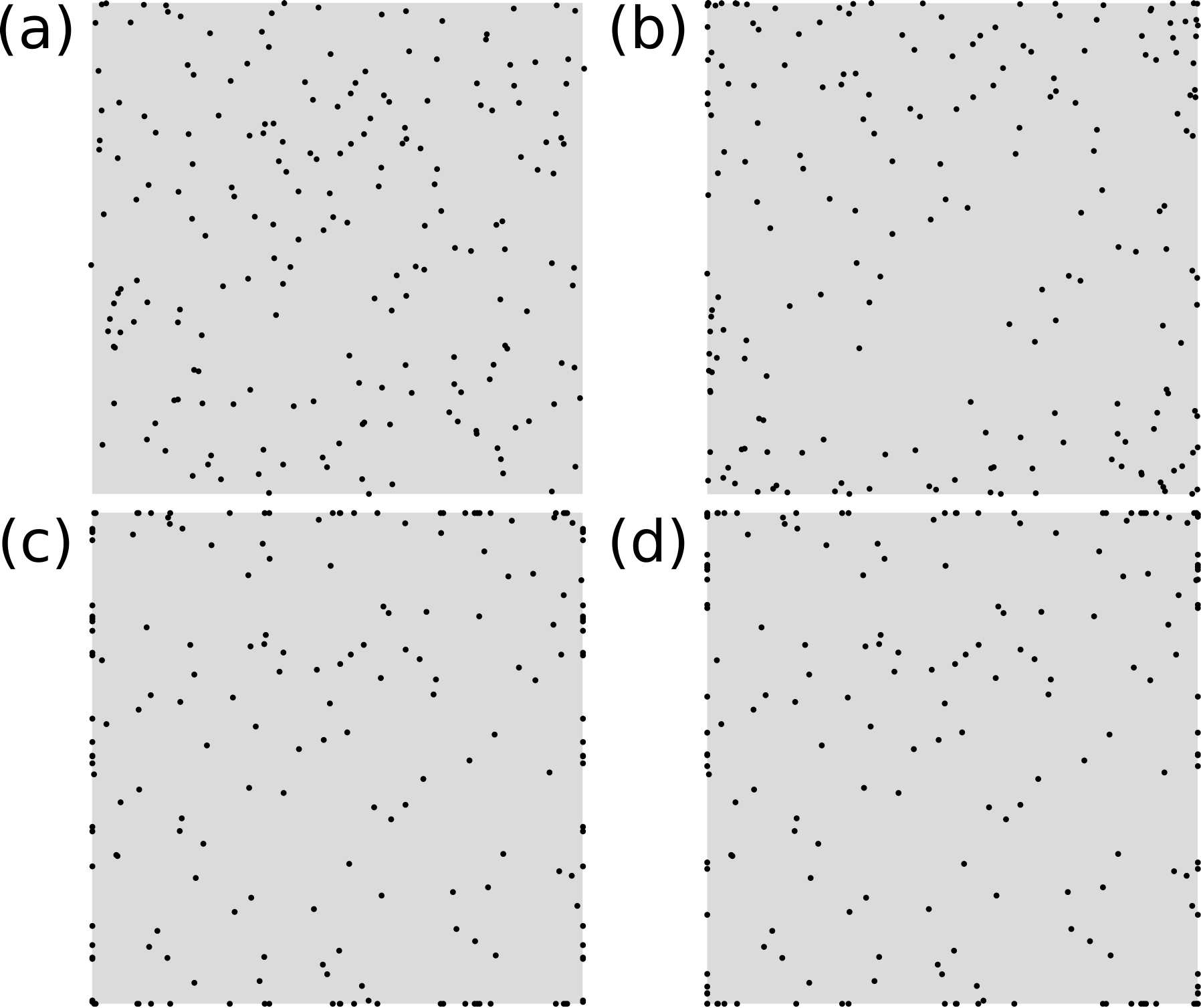}
\caption{Examples of sampling distributions on the square: (a) $\nu$, (b) $\nu'$,
(c) $\nu_{1/2}$, (d) $\nu'_{1/2}$.}
\label{squarespl}
\end{figure}

We numerically compute $K(m)$ for four different distributions:
\begin{itemize}
\item  $\nu_0 = ds$, the uniform distribution on the square,
\item  $\nu' = 4ds/\left(\pi^2 \sqrt{1 - x^2}\sqrt{1 - y^2}\right)$, a distribution denser near the edges and corners of the square,
\item $\nu_\alpha = (1-\alpha)ds + \alpha d\sigma$, where $\sigma$ is the uniform
distribution on the boundary of the square,
\item $\nu'_\alpha = (1-\alpha)ds + \alpha d\sigma'$, where $d\sigma'$ is the measure on
the boundary with weight  $1/4\pi\sqrt{1 - s^2}$ where $s = \min (x,y)$ (i.e.\ denser near the corners of the square).
\end{itemize}
Figure \ref{squarespl} shows examples of such  distributions for $\alpha = 1/2$.

\begin{figure}
\centering
\includegraphics[width=10cm]{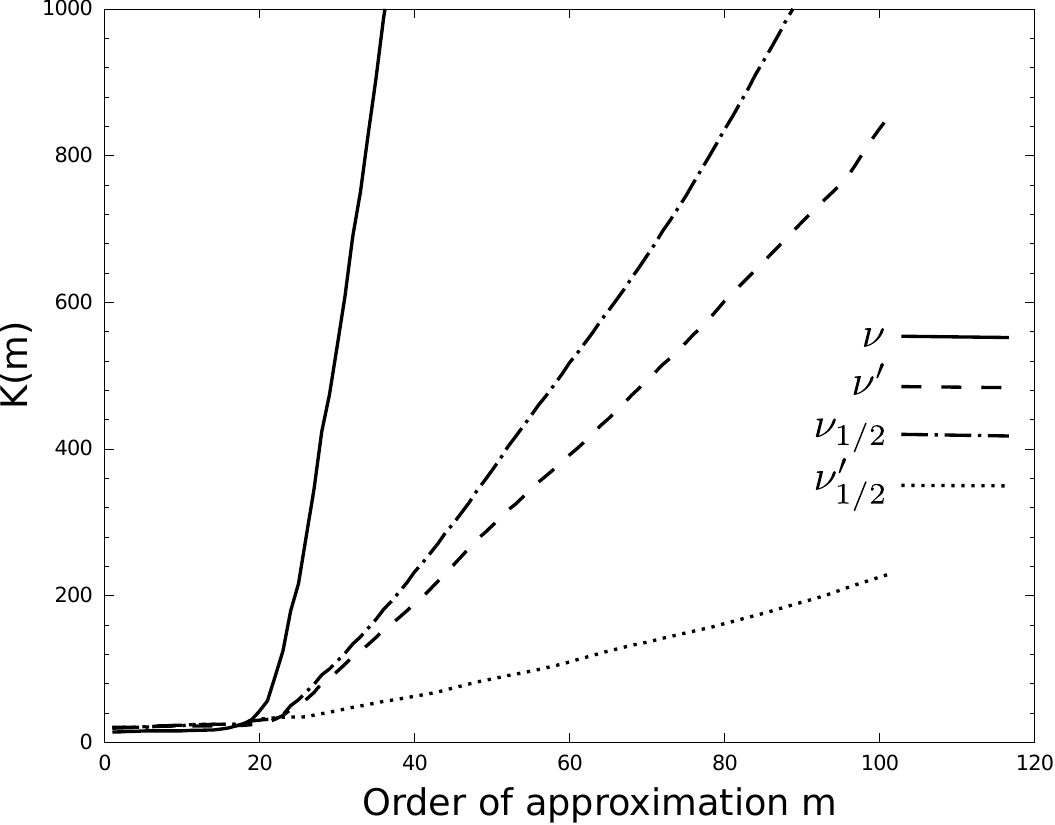}
\caption{Numerical evaluation of $K$ for four different samples distribution
on the square.}
\label{square}
\end{figure}

The estimated values of $K(m)$ for
 $\nu_0$, $\nu'$, $\nu_{1/2}$ and  $\nu'_{1/2}$ are given on figure~\ref{square}. Here, having a denser sampling near or on the border of the square improves the stability of the reconstruction
compared to the uniform case, but still needs a high number of samples.

Using the non-uniform sampling on the border, with more samples in
the sections of the boundary furthest from the origin, makes the behavior
of $K(m)$ comparable to $m$, which is the best case possible.

\section{Conclusion}

In this paper, we  compare different ways of sampling solutions to the
Helmholtz equation, using a finite number of point measurements.
Our main results reveal that good reconstructions
can be obtained using Fourier-Bessel or
plane waves approximations, and that these
reconstructions benefit from a denser sampling 
on the boundary.

These results were obtained in the particular case of a two-dimensional disc
or a three-dimensional ball.
For a more general star-shaped domain in $\R^2$, the Fourier-Bessel
approximation remains valid, but the quantity $K(m)$ does not
have an explicit expression, yet it can be evaluated
numerically after orthogonalization of the Fourier-Bessel or plane waves family.
Our first numerical investigation, on the square, indicates that denser, but non-uniform, sampling of the functions
on the boundary is also beneficial in this more general setting.

Finally, sampling of other physical quantities can benefit from similar sampling strategies,
e.g.\ vibrations of plates \cite{CLD, CD}, electromagnetic fields \cite{HMP3}, or vibrations in 3D linear elasticity \cite{Mo}.
Indeed, they can be approximated using schemes similar to the ones used here.

\vspace{13pt}
\centerline{ACKNOWLEDGEMENT}
\vspace{13pt}

 This research is supported by the FWF START-project FLAME (Y 551-N13) and
the ANR project ECHANGE (ANR-08-EMER-006). LD is partly supported by Institut Universitaire de France and by LABEX WIFI (Laboratory of
Excellence within the French Program "Investments for the Future") under
references ANR-10-LABX-24 and ANR-10-IDEX-0001-02 PSL.


\end{document}